\documentclass[12pt]{article}
\usepackage{amssymb}
\title{On weak approximation of $U$-statistics}
\author{Masoud M. Nasari\footnote{\footnotesize{Research supported by a Carleton university Faculty of Graduate Studies and Research scholarship, and NSERC Canada Discovery Grants of M. Cs\"{o}rg\H{o} and M. Mojirsheibani at Carleton university.}} \\ \footnotesize	{\emph{School of Mathematics and Statistics, Carleton University, Canada}}\\
\footnotesize{\emph{e-mail}: \texttt{mmnasari@connect.carleton.ca}}}
\date{}
\begin{document}
\maketitle
\abstract{$ $ \\ This paper investigates weak convergence of $U$-statistics via approximation in probability. The classical condition that the second moment of the kernel of the underlying $U$-statistic exists is relaxed to having $\frac{4}{3}$ moments
  only (modulo a logarithmic term). Furthermore, the conditional
  expectation of the kernel is only assumed to be in the domain of
  attraction of the normal law (instead of
the classical two-moment condition).\\}
\section{Introduction}
$ $ \\
 Employing truncation arguments and the concept of weak convergence of self-normalized and studentized partial sums, which were inspired by the works of
Cs\"{o}rg\H{o}, Szyszkowicz and Wang in \cite{cs1}, \cite{cs5}, \cite{ana} and \cite{cs}, we derive weak convergence
   results via approximations in probability for pseudo-self-normalized U-statistics and U-statistic type  processes. Our results require only that (i) the expected value of the product of the kernel of the underlying $U$-statistic to the exponent $\frac{4}{3}$ and its logarithm exists
(instead of having 2 moments of the kernel), and  that (ii) the conditional expected value of the kernel on each observation is in the domain of attraction of the normal law (instead of having 2 moments). Similarly relaxed moment conditions were first used by Cs\"{o}rg\H{o}, Szyszkowicz and Wang \cite{cs1} for $U$-statistics type processes  for changepoint problems in terms of kernels of order 2 (cf. Remark 5). Our  results in this exposition extend their work to approximating $U$-statistics with higher order kernels. The thus obtained  weak convergence results for $U$-statistics in turn extend those  obtained by R.G. Miller Jr. and P.K. Sen in \cite{mil} in 1972 (cf. Remark 3). The latter results of Miller and Sen are based on the classical condition of  the existence of the second moment of the kernel of the underlying $U$-statistic which in turns implies the existence of the second moment of the conditional expected value of the kernel on each of the  observations.
\newpage
\section{ Main results and Background}
Let $X_{1},X_{2},\ldots$, be a sequence of non-degenerate real-valued i.i.d. random variables with distribution $F$. Let
$h(X_{1},\ldots,X_{m})$, symmetric in its arguments, be a Borel-measurable real-valued kernel of order $m\geq 1$, and consider the
parameter $\theta=\underbrace{\int\ldots\int}_{\mathbb{R}^m} h(x_{1},\ldots ,x_{m})\ dF(x_{1})\ldots dF(x_{m})<\infty.$ The corresponding $U$-statistic (cf. Serfling  \cite{serf}
or Hoeffding \cite{hof}) is
$$U_{n}={n \choose m}^{-1}
\sum_{C(n,m)}\ h(X_{i_{1}},\ldots,X_{i_{m}}),$$
 where $m\leq n$ and $\sum_{C(n,m)}$ denotes the sum over $C(n,m)=\{1\leq i_{1}<\ldots <i_{m}\leq n\}$.
\\
\par
In order to state our results, we first need the following definition.\\
\textbf{Definition}. A sequence $X,X_{1},X_{2},\ldots,$ of i.i.d. random variables is said to be in the domain of attraction of the normal law ($X\in$ $DAN$) if there exist sequences of constants $A_{n}$ and $B_{n}>0$ such that, as $n\rightarrow \infty,$
$$\frac{\sum_{i=1}^{n}X_{i}-A_{n}}{B_{n}}\longrightarrow_{d} N(0,1).$$
\textbf{Remark 1}. Furtherer to this definition of $DAN$, it is known that $A_{n}$ can be taken as $n\mathbb{E}(X)$ and $B_{n}=n^{1/2} \ell_{X}(n)$, where $\ell_{X}(n)$ is a slowly varying function at infinity (i.e., $\lim_{n\rightarrow \infty}
\frac{\ell_{X}(nk)}{\ell_{X}(n)}=1$ for any $k>0$), defined by the distribution of $X$. Moreover,
$\ell_{X}(n)=\sqrt{Var(X)}>0$, if $Var(X)< \infty$, and $\ell_{X}(n)\rightarrow \infty$, as $n\rightarrow \infty$, if $Var(X)=\infty$. Also $X$ has all moments less than 2, and the variance of $X$ is positive, but need not be finite.
\newline \\
Also define the pseudo-self-normalized $U$-process as follows.\\ $$U_{[nt]}^{*}=\left\{
\begin{array}{ll}
\ 0 \qquad \qquad  \  ,& \hbox{$0\leq t <\displaystyle{\frac{m}{n}}, $}\\
                     \displaystyle{\frac{U_{[nt]}-\theta }{V_{n}}} \  \ \ \  , & \hbox{$\displaystyle{\frac{m}{n}}\leq                      t\leq 1,$}\\
                    \end{array}
                  \right.$$\\ where [.] denotes the greatest integer function,  $V_{n}^{2}:=\sum_{i=1}^{n}\tilde{h}^{2}_{1}(X_{i})$ and $\tilde{h}_{1}(x)=\mathbb{E}\textbf{(}h(X_{1},\ldots,X_{m})-\theta | X_{1} =x\textbf{)}$.
\\ \\
\textbf{Theorem 1.} $If$ \\ \\
(a) $\mathbb{E}\left(|h(X_{1},\dots,X_{m})|^{\frac{4}{3}}\log|h(X_{1},\ldots,X_{m})|\right)<\infty \ and$
$\tilde{h}_{1} (X_{1}) \in DAN$,\\

$then,\ as\ n\rightarrow \infty,\ we \ have$\\ \\
(b) $\displaystyle{\frac{[nt_{0}]}{m}\ U^{*}_{[nt_{0}]}}\rightarrow_{d} \ N(0,t_{0}), \ for \ t_{0}\in(0,1]$;\\ \\
(c) $\displaystyle{\frac{[nt]}{m}\ U^{*}_{[nt]}}\ \rightarrow_{d}\ W(t)$ \emph{on} ($D$[0,1],$\rho$), $where$ $\rho$ \emph{is the sup-norm for functions in}\\

 $D[0,1] \ and\  \{W(t), 0\leq t \leq 1\}\ is\ a\ standard\ Wiener\ process$;\\ \\
(d) $On\ an\ appropriate\ probability \ space \ for\ X_{1},X_{2},\ldots,\ we\ can\ construct\ a$

$standard \ Wiener \ process\ \{W(t), 0\leq t < \infty\}\ such\ that$
$$\sup_{0\leq t \leq 1}\ \left|\ \frac{[nt]}{m}\ U^{*}_{[nt]} -\ \frac{W(nt)}{n^{\frac{1}{2}}} \right|=o_{P}(1).$$ 
\textbf{Remark 2}. The statement (c), whose notion will be used throughout, stands for the following functional central limit theorem (cf. Remark 2.1 in Cs\"{o}rg\H{o}, Szyszkowicz and Wang  \cite{cs}). On account of (d), as $n\rightarrow\infty$, we have
$$g(S_{[n.]}/V_{n})\longrightarrow_{d} g(W(.))$$
for all $g:D=D[0,1]\longrightarrow \mathbb{R}$ that are $(D,\mathfrak{D})$ measurable and  $\rho$-continuous, or $\rho$-continuous except at points forming a set of Wiener measure zero on $(D,\mathfrak{D})$, where $\mathfrak{D}$ denotes the $\sigma$-field of subsets of $D$ generated by the finite-dimensional subsets of $D$.\\ 
\\
\par
Theorem 1 is fashioned after the work on  weak convergence of  self-normalized partial sums processes of Cs\"{o}rg\H{o}, Szyszkowicz and Wang in \cite{ana}, \cite{cs} and \cite{cs5}, which constitute extensions of the contribution of Gin\'{e}, G\"{o}tze and Mason in \cite{Gin}. \\ \\
\par
As to $\tilde{h}_{1}(X_{1}) \in DAN$, since $\mathbb{E} \tilde{h}_{1}(X_{1})=0$ and $\tilde{h}_{1}(X_{1}),\tilde{h}_{1}(X_{2}),\ldots,$ are i.i.d. random variables, Theorem 1 of \cite{ana} (cf. also Theorem 2.3 of \cite{cs}) in this context reads as follows.
\\
\textbf{Lemma 1.} \emph{As $n\rightarrow \infty,$ the following statements are equivalent}:
\\ \\
(a) \ \ $\tilde{h}_{1}(X_{1})\in$ $DAN$ ;\\ \\
(b) $\ \displaystyle{\frac{\sum_{i=1}^{[nt_{0}]}\tilde{h}_{1}(X_{i})}{V_{n}}}\longrightarrow_{d} N(0,t_{0})$ \emph{for} $t_{0}\in(0,1];$ \\\\ \\
$(\textrm{c}) \ \displaystyle{\frac{\sum_{i=1}^{[nt]}\tilde{h}_{1}(X_{i})}{V_{n}}}\longrightarrow_{d} W(t)$ \emph{on} $(D[0,1],\rho)$,
\emph{where $\rho$ is the sup-norm metric}

\emph{for functions in} $D[0,1]$ \emph{and} $\{W(t),0\leq t\leq 1\}$ \emph{is a standard Wiener}

\emph{process};\\\\
(d) \emph{On an appropriate probability space for} $X_{1},X_{2},\ldots,$ \emph{we can construct a}

\emph{standard Wiener process} $\{W(t),0\leq t<\infty\}$ \emph{such that}
$$\sup_{0\leq t \leq 1}\left|\frac{\sum_{i=1}^{[nt]}\tilde{h}_{1}(X_{i})}{V_{n}} -\ \frac{W(nt)}{n^{\frac{1}{2}}} \right|=o_{P}(1).$$
\\ \\
Also, in the same vein, Proposition 2.1 of \cite{cs} for $\tilde{h}_{1}(X_{1})\in$ DAN reads as follows.\\ \\
\textbf{Lemma 2.} \emph{As} $n\rightarrow \infty$,  \emph{the following statements are equivalent}: \\ \\
(a)\ \ $\tilde{h}_{1}(X_{1})\in$ $DAN$; \\

\emph{There is a sequence of constants} $B_{n}\nearrow \infty$, \emph{such that}\\ \\
(b) $\ \displaystyle{\frac{\sum_{i=1}^{[nt_{0}]}\tilde{h}_{1}(X_{i})}{B_{n}}}\longrightarrow_{d} N(0,t_{0})$ \emph{for} $t_{0}\in(0,1];$ \\\\ \\
$(\textrm{c}) \ \displaystyle{\frac{\sum_{i=1}^{[nt]}\tilde{h}_{1}(X_{i})}{B_{n}}}\longrightarrow_{d} W(t)$ \emph{on} $(D[0,1],\rho)$,
\emph{where} $\rho$ \emph{is the sup-norm metric}

\emph{for functions in} $D[0,1]$ \emph{and} $\{W(t),0\leq t\leq 1\}$
\emph{is a standard Wiener}

\emph{process};\\\\
(d) \emph{On an appropriate probability space for} $X_{1},X_{2},\ldots,$ \emph{we can construct a}

\emph{standard Wiener process} $\{W(t),0\leq t<\infty\}$ \emph{such that}\\
$$\sup_{0\leq t \leq 1}\left|\frac{\sum_{i=1}^{[nt]}\tilde{h}_{1}(X_{i})}{B_{n}} -\ \frac{W(nt)}{n^{\frac{1}{2}}} \right|=o_{P}(1).$$
\\
In view of Lemma 2, a scalar normalized companion of Theorem 1 reads as follows.\\ \\
\textbf{Theorem 2.} \emph{If}
\\ \\
(a) $\mathbb{E}\left(|h(X_{1},\dots,X_{m})|^{\frac{4}{3}}\log|h(X_{1},\ldots,X_{m})|\right)<\infty$  \emph{and}
$\tilde{h}_{1} (X_{1}) \in$ $DAN$,\\

\emph{then}, as $n\rightarrow \infty$, \emph{we have}\\ \\
(b) $\displaystyle{\frac{[nt_{0}]}{m}\ \frac{U_{[nt_{0}]}-\theta }{B_{n}}}\longrightarrow_{d} \ N(0,t_{0}),\ where \ t_{0}\in(0,1]$;\\ \\
(c) $\displaystyle{\frac{[nt]}{m}\ \frac{U_{[nt]}-\theta }{B_{n}}}\ \longrightarrow_{d}\ W(t)$ \emph{on} ($D$[0,1],$\rho$), \emph{where} $\rho$ \emph{is the sup-norm for} \\

\emph{functions in}   $D$[0,1] \emph{and} $\{W(t), 0\leq t \leq 1\}$ \emph{is a standard Wiener process};\\ \\
(d) \emph{On an appropriate probability space for} $X_{1},X_{2},\ldots$, \emph{we can construct a}

\emph{standard Wiener process} $\{W(t), 0\leq t < \infty\}$ \emph{such that}
$$\sup_{0\leq t \leq 1}\ \left|\ \frac{[nt]}{m}\frac{U_{[nt]}-\theta }{B_{n}} -\ \frac{W(nt)}{n^{\frac{1}{2}}} \right|=o_{P}(1).$$
\\\\\\\\
By defining
\begin{eqnarray*}
Y_{n}^{*}(t)&=&0\ \ \qquad  \qquad \qquad \qquad  \qquad \qquad \qquad  \textrm{for} \
0\leq t\leq \frac{m-1}{n},\\ Y^{*}_{n}(\frac{k}{n})&=&\frac{k(U_{k}-\theta)}{m \sqrt{n
Var\textbf{(}\tilde{h}_{1}(X_{1})\textbf{)}}} \qquad \qquad \qquad
\ \ \textrm{for}\ k=m,\ldots,n
\end{eqnarray*}
and for $t\in\left[\frac{k-1}{n}\ ,\ \frac{k}{n}\right],\  k=m,\ldots,n$ ,
$$Y^{*}_{n}(t)\ =\ Y^{*}_{n}(\frac{k-1}{n})+n(t-\frac{k-1}{n})\left(Y^{*}_{n}(\frac{k}{n})-Y^{*}_{n}(\frac{k-1}{n})\right),\ \ \  \ $$
we can state the already mentioned 1972 weak convergence result of Miller and Sen as follows.\\\\
\textbf{Theorem A}.
\emph{If} $$(\textrm{I})\  0<\mathbb{E}\textbf{[}\textbf{(}h(X_{1},X_{2},\ldots ,X_{m})-\theta\textbf{)}\textbf{(}h(X_{1},X_{m+1},\ldots
,X_{2m-1})-\theta\textbf{)}\textbf{]}=Var(\tilde{h}_{1}(X_{1}))< \infty$$

\emph{and}
$$(\textrm{II}) \ \mathbb{E}h^{2}(X_{1},\ldots,X_{m})< \infty,\qquad \qquad \qquad \qquad \qquad \qquad \qquad \qquad \qquad \qquad \qquad \qquad \qquad$$

\emph{then, as} $n\rightarrow \infty,$\\
$$Y_{n}^{*}(t)\rightarrow_{d} W(t)\ \ \ \textrm{on}\  (\textrm{C}[0,1],\rho),$$
\emph{where} $\rho$ \emph{is the sup-norm for functions in} $\textrm{C}[0,1]$ \emph{and} $\{W(t),0\leq t \leq 1\}$ \emph{is a standard Wiener process} .\\
\newline
\\
\textbf{Remark 3.} When $\mathbb{E}h^{2}(X_{1},\ldots,X_{m})< \infty$, first note that existence of  the second moment of the kernel $h(X_{1},\ldots,X_{m})$ implies the existence of the second moment of $\tilde{h}_{1}(X_{1})$. Therefore, according to Remark 1, $B_{n}=\sqrt{n\ \mathbb{E}\tilde{h}_{1}^{2}(X_{1})}$. This means that under the  conditions of Theorem A, Theorem 2 holds true and, via (c) of latter, it yields a version of Theorem A on $D[0,1]$. We note in passing that our method of proofs differs from that of cited paper of Miller and Sen. We use a method of truncation \`{a} la \cite{cs1} to relax the condition $\mathbb{E}h^{2}(X_{1},\ldots,X_{m})<\infty$ to the less stringent moment condition $\mathbb{E}\left(|h(X_{1},\dots,X_{m})|^{\frac{4}{3}}\log|h(X_{1},\ldots,X_{m})|\right)<\infty$ that, in turn, enables us to have $\tilde{h}_{1}(X_{1})\in DAN$ in general, with the possibility of infinite variance.
\\ \\
\textbf{Remark 4.} Theorem 1 of \cite{ana} (Theorem 2.3 in \cite{cs}) as well as Proposition 2.1 of \cite{cs}, continue to hold true in terms of Donskerized partial sums that are elements of $C[0,1].$ Consequently, the same is true for the above stated Lemmas 1 and 2, concerning $\tilde{h}_{1}(X_{1})\in DAN$. This in turn, mutatis mutandis, renders appropriate versions of Theorems 1 and 2 to hold true in $(C[0,1],\rho).$\\\\\\
\textbf{Proof of Theorems 1 and 2.} \\
In view of Lemmas 1 and 2, in order to prove Theorems 1 and 2, we only have to prove the following theorem.\\ \\
\textbf{Theorem 3.} \emph{If} $\mathbb{E}\left(|h(X_{1},\dots,X_{m})|^{\frac{4}{3}}\log|h(X_{1},\ldots,X_{m})|\right)<\infty$ \emph{and}  $\tilde{h}_{1} (X_{1})\in$ \emph{DAN then, as} $n\rightarrow \infty$, \emph{we have}\\\\
  $ \ \sup_{0\leq t \leq 1} \left|\displaystyle{\frac{[nt]}{m}} U^{*}_{[nt]}- \displaystyle{\frac{\sum_{i=1}^{[nt]}\tilde{h}_{1}(X_{i})}{V_{n}}} \right|=o_{P}(1),$\qquad \qquad \qquad \qquad \qquad\qquad \qquad  \textbf{(1)}\\
\emph{and}\\
 $\sup_{0\leq t \leq 1} \left|\displaystyle{\frac{[nt]}{m}} \frac{U_{[nt]}-\theta}{B_{n}}- \displaystyle{\frac{\sum_{i=1}^{[nt]}\tilde{h}_{1}(X_{i})}{B_{n}}} \right|=o_{P}(1).$\qquad \qquad \qquad \qquad \qquad \qquad  \textbf{(2)}
\\\\
\textbf{Proof of Theorem 3}. In view of  (b) of Lemma 2 with $t_{0}=1$, Corollary 2.1 of \cite{cs}, yields $\displaystyle{\frac{V_{n}^{2}}{B_{n}^{2}}}\rightarrow_{P} 1$. This in turn implies the equivalency of (1) and (2). Therefore, it suffices to prove (2) only.
\\ \\
It can be easily seen that\\
$$\sup _{0\leq t\leq 1} \left|\ \displaystyle{\frac{[nt]}{m} \frac{U_{[nt]}-\theta}{B_{n}}} - \frac{\sum_{i=1}^{[nt]} \tilde{h}_{1}(X_{i})}{B_{n}}\right|\leq \sup_{0\leq t < \frac{m}{n}} \left| \frac{\sum_{i=1}^{[nt]} \tilde{h}_{1}(X_{i})}{B_{n}}\right|\qquad \qquad \qquad \qquad$$ \\
$$\qquad \qquad \qquad \qquad \qquad  \qquad \qquad \qquad \qquad  +\sup_{\frac{m}{n} \leq t \leq 1} \left| \frac{[nt]}{m}\frac{U_{[nt]}-\theta}{B_{n}}- \frac{\sum_{i=1}^{[nt]} \tilde{h}_{1}(X_{i})}{B_{n}}\right|.$$\\
 Since, as $n\rightarrow\infty$, we have $\displaystyle{ \frac{m}{n}}\rightarrow 0$ and, consequently, in view of (d) of Lemma 2\\
$$\sup_{0\leq t < \frac{m}{n}}\left| \displaystyle{\frac{\sum_{i=1}^{[nt]} \tilde{h}_{1}(X_{i})}{B_{n}}}\right|=o_{P}(1),$$ in order to prove (2), it will be enough to show that
$$\sup_{\frac{m}{n} \leq t \leq 1} \left| \frac{[nt]}{m}\frac{U_{[nt]}-\theta}{B_{n}}- \frac{\sum_{i=1}^{[nt]} \tilde{h}_{1}(X_{i})}{B_{n}}\right|=o_{P}(1),$$
or equivalently to show that
 \\$$\max_{m\leq k \leq n}\left|\ \frac{k}{m B_{n}}{k \choose m}^{-1} \sum_{C(k,m)} \textbf{(}h(X_{i_{1}},\ldots,X_{i_{m}})-\theta\textbf{)}-\frac{1}{B_{n}} \sum_{i=1}^{k}\tilde{h}_{1}(X_{i})\ \right|\qquad \qquad \qquad$$
\begin{eqnarray*}
&=&\max_{m\leq k \leq n}\left | \ \frac{k}{m B_{n}}{k \choose m}^{-1} \sum_{C(k,m)}\
                                                                                                  \left(h(X_{i_{1}},\ldots,X_{i_{m}})-\theta-\tilde{h}_{1}(X_{i_{1}})-\ldots-\tilde{h}_{1}(X_{i_{m}})\right)\ \right | \\
\\ &=& o_{P}(1).\qquad \qquad \qquad \qquad \qquad \qquad \qquad \qquad \qquad \qquad \qquad \qquad
\qquad \qquad \qquad  \textbf{(3)}
\end{eqnarray*}
The first equation of (3) follows from the fact that \\
$$\sum_{C(k,m)} \left(\tilde{h}_{1}(X_{i_{1}})+ \ldots +\tilde{h}_{1}(X_{i_{m}})\right)=\frac{m}{k}  {k \choose m} \sum_{i=1}^{k} \tilde{h}_{1}(X_{i}),$$
\\ where $\sum_{C(k,m)}$ denotes the sum over $C(k,m)=\{1\leq i_{1} < \ldots < i_{m} \leq k \}$. To establish (3), without loss of generality we can, and shall assume that $\theta=0$.\\ \\
Considering that for large $n$, $\displaystyle{\frac{1}{B_{n}}}\leq \displaystyle{\frac{1}{\sqrt{n}}}$ (cf. Remark 1), to conclude (3), it will be enough to show that, as $n\rightarrow \infty$, the following holds: \\
$$n^\frac{-1}{2} \max _{m \leq k \leq n}  \left|k {k \choose m}^{-1}
 \sum _{C(k,m)} \left( h(X_{i_{1}}, \ldots ,X_{i_{m}})-\tilde{h}_{1}(X_{i_{1}})- \ldots -\tilde{h}_{1}(X_{i_{m}})\right) \right|
=o_{P}(1). \ \  \textbf{(4)}$$
\par
To establish (4), for the ease of notation, let \\
$$h^{(1)}(X_{i_{1}},
\ldots, X_{i_{m}}) := h(X_{i_{1}},\ldots,X_{i_{m}}) I_{(|h|\leq n^{\frac{3}{2}})}-\mathbb{E}\textbf{(}h(X_{i_{1}},\ldots,X_{i_{m}}) I_{(|h|\leq n^{\frac{3}{2}})}\textbf{)},$$
$$\tilde{h}^{(1)}(X_{i_{j}}):=\mathbb{E}\textbf{(}h^{(1)}(X_{i_{1}}, \ldots, X_{i_{m}})|X_{i_{j}}\textbf{)},   \ \ \ j=1,\ldots,m,\qquad \qquad \qquad \qquad \qquad \qquad $$
$$\psi^{(1)}(X_{i_{1}},\ldots,X_{i_{m}}):=h^{(1)}(X_{i_{1}},\ldots,X_{i_{m}})-\tilde{h}^{(1)}(X_{i_{1}})-\ldots-\tilde{h}^{(1)}(X_{i_{m}}),\qquad \qquad \qquad \ $$
$$h^{(2)}(X_{i_{1}}, \ldots, X_{i_{m}})
  := h(X_{i_{1}}, \ldots, X_{i_{m}})I_{(|h|> n^{\frac{3}{2}})}-\mathbb{E}\textbf{(}h(X_{i_{1}}, \ldots, X_{i_{m}})I_{(|h|> n^{\frac{3}{2}})}\textbf{)},\ \ $$
$$\tilde{h}^{(2)}(X_{i_{j}})
 :=\mathbb{E}\textbf{(}h^{(2)}(X_{i_{1}}, \ldots, X_{i_{m}})|X_{i_{j}}\textbf{)}, \  j=1,\ldots,m,\qquad \qquad \qquad \qquad \qquad \qquad \ \ \ \qquad $$\\
where $I_{A}$ is the indicator function of the set $A$.
Now observe that\\
\\ $$n^\frac{-1}{2} \max _{m \leq k \leq n}  \left| k {k \choose m}^{-1}
\sum _{C(k,m)} \left(h(X_{i_{1}},\ldots,X_{i_{m}})-\tilde{h}_{1}(X_{i_{1}})- \ldots -\tilde{h}_{1}(X_{i_{m}})\right) \right|\ \ \ \ \ \ \qquad \ $$\\
$$\leq \ \  n^\frac{-1}{2} \max _{m \leq k \leq n}  \left| k {k \choose m}^{-1}
\sum _{C(k,m)} \left( h(X_{i_{1}},\ldots,X_{i_{m}})-h^{(1)}(X_{i_{1}},\ldots,X_{i_{m}})\right) \right|\qquad \qquad \qquad$$
$$\ \  +\  n ^\frac{-1}{2} \max _{m \leq k \leq n}  \left| k {k \choose m}^{-1}
\sum _{C(k,m)} \left( \tilde{h}_{1}(X_{i_{1}})+ \ldots + \tilde{h}_{1}(X_{i_{m}})-
\tilde{h}^{(1)}(X_{i_{1}})- \ldots -
\tilde{h}^{(1)}(X_{i_{m}}) \right) \right| $$
$$\ \ + \   n ^\frac{-1}{2} \max _{m \leq k \leq n}  \left| k {k \choose m}^{-1}
\sum _{C(k,m)}\ \psi^{(1)}(X_{i_{1}},\ldots,X_{i_{m}})\ \right|\ \ \ \ \ \ \ \ \ \ \ \ \ \ \ \ \ \ \ \ \ \ \ \qquad \qquad \  \ \ \ \ \ \  $$ \\
$ $ $\ \ \ \ \ :=J_{1}(n)+J_{2}(n)+J_{3}(n).$ \\ \\
We will show that $J_{s}(n)=o_{P}(1)$, $s=1,2,3.$\\
\par
To deal with the term $J_{1}(n)$, first note that
\\
$$ h(X_{i_{1}},\ldots,X_{i_{m}})-h^{(1)}(X_{i_{1}},\ldots,X_{i_{m}})=h^{(2)}(X_{i_{1}},\ldots,X_{i_{m}}).$$ \\
Therefore, in view of Theorem 2.3.3 of \cite{boro} page 43, for $\epsilon>0$,  we can write
\\
$$\mathbb{P}\left( n^\frac{-1}{2} \max _{m \leq k \leq n}  \left| k {k \choose m}^{-1}
\sum _{C(k,m)} h^{(2)}(X_{i_{1}},\ldots,X_{i_{m}}) \right| > \epsilon\right)\ \ \ \qquad \qquad \qquad \qquad \ \ \ \ \   $$ \\
$\leq \epsilon^{-1} n^{\frac{-1}{2}} \left(\ m\ \mathbb{E}| h^{(2)}(X_{1},\ldots,X_{m})|+ \ n\ \mathbb{E}| h^{(2)}(X_{1},\ldots,X_{m})|\ \right)$\\ \\
$\leq \epsilon^{-1} n^{\frac{-1}{2}} \ 2m \ \mathbb{E}|h(X_{1},\ldots,X_{m})| + \epsilon^{-1} n^{\frac{1}{2}}\ 2m\  \mathbb{E}\textbf{(}|h(X_{1},\ldots,X_{m})| I_{(|h|> n^{\frac{3}{2}})}\textbf{)}$\\
\\
$\leq \epsilon^{-1} n^{\frac{-1}{2}} \ 2m\  \mathbb{E}|h(X_{1},\ldots,X_{m})| + \epsilon^{-1} \ 2m \ \mathbb{E} \textbf{(}|h(X_{1},\ldots,X_{m})|^{\frac{4}{3}} I_{( |h|> n^{\frac{3}{2}})}\textbf{)}$
$$\longrightarrow 0, \ \ \ \ \textrm{as} \ n\rightarrow \infty.\qquad \qquad \qquad \qquad \qquad \qquad \qquad \qquad \qquad \qquad \qquad \qquad \qquad $$
Here we have used the fact that $\mathbb{E}|h(X_{1},\ldots,X_{m})|^{\frac{4}{3}} < \infty$. The last line above implies that $J_{1}(n)=o_{P}(1).$
\\
\par
Next to deal with $J_{2}(n)$, first observe that \\
$$\tilde{h}_{1}(X_{i_{1}})+ \ldots
+\tilde{h}_{1}(X_{i_{m}})-\tilde{h}^{(1)}(X_{i_{1}})- \ldots -
\tilde{h}^{(1)}(X_{i_{m}}) =\sum_{j=1}^{m} \tilde{h}^{(2)}(X_{i_{j}}).$$ \\
It can be easily seen that $\sum_{j=1}^{m} \tilde{h}^{(2)}(X_{i_{j}})$ is symmetric in $X_{i_{1}}, \ldots ,X_{i_{m}}$. Thus, in view of Theorem 2.3.3 of \cite{boro} page 43, for $\epsilon>0$, we have
$$\mathbb{P}\left( n^{\frac{-1}{2}} \max_{m\leq k \leq n}k \left| {k \choose m}^{-1} \sum_{C(k,m)}
\left(\sum_{j=1}^{m} \tilde{h}^{(2)}(X_{i_{j}})\right) \right|>\epsilon\right)\qquad \qquad \qquad \qquad \qquad \qquad \qquad \qquad$$
$\leq \epsilon^{-1} n^{\frac{-1}{2}}\ 2m\ \mathbb{E}|h(X_{1},\ldots,X_{m})| +\epsilon^{-1} n^{\frac{1}{2}}\  2m\   \mathbb{E}\textbf{(}|h(X_{1},\ldots,X_{m})|I_{(|h|> n^{\frac{3}{2}})}\textbf{)}$
$$\longrightarrow 0, \ \ \ \ \textrm{as} \ n\rightarrow \infty,\qquad \qquad \qquad \qquad \qquad \qquad \qquad \qquad \qquad \qquad \qquad \qquad \qquad $$  i.e., $J_{2}(n)=o_{P}(1).$ \\ \\
\textbf{Note}. Alternatively, one can use Etemadi's maximal inequality for partial sums of i.i.d. random variables, followed by Markov inequality, to show $J_{2}(n)=o_{P}(1).$
\\ \\ \\\\\\
\par
As for the term $J_{3}(n)$, first note that  ${k \choose m}^{-1}
\sum_{C(k,m)} \psi^{(1)}(X_{i_{1}},\ldots,X_{i_{m}})
$ is a $U$-statistic. Consequently one more application of Theorem 2.3.3 page 43 of \cite{boro} yields,
$$\mathbb{P}\left(n^{\frac{-1}{2}} \max_{m \leq k \leq n} k \left|  {k \choose m}^{-1}
\sum _{C(k,m)} \psi^{(1)}(X_{i_{1}},\ldots,X_{i_{m}}) \right|>\epsilon\right)\qquad \qquad \qquad \qquad \qquad \qquad \qquad \qquad$$
$\leq n^{-1} \epsilon^{-2} \  m^{2} \ \mathbb{E}\textbf{(}\psi^{(1)}(X_{1},\ldots,X_{m})\textbf{)}
^{2} $ \\
$$ \ \ \ +\  n^{-1} \epsilon^{-2} \ \sum_{k=m+1}^{n}(2k+1)\   \mathbb{E}\left( {k \choose m}^{-1}
\sum _{C(k,m)} \psi^{(1)}(X_{i_{1}},\ldots,X_{i_{m}})\right)^{\displaystyle{2}}.\ \ \ \qquad \ \  \  \textbf{(5)}$$\\
Observing that $\mathbb{E}\textbf{(}\psi^{(1)}(X_{1},\ldots,X_{m})\textbf{)}^{2}\leq C(m) \ \mathbb{E}\left(h^{2}(X_{1},\ldots,X_{m}) I_{(|h|\leq n^{\frac{3}{2}})}\right),$ where $C(m) $ is a positive constant that does not depend on $n$,  $$\mathbb{E}\psi^{(1)}(X_{i_{1}},\ldots,X_{i_{m}})=\mathbb{E}\textbf{(}\psi^{(1)}(X_{i_{1}},\ldots,X_{i_{m}})|X_{i_{j}}\textbf{)}=0,\ j=1,\dots,m,$$
and in view of Lemma B page 184 of \cite{serf}, it follows that for some positive constants $C_{1}(m)$ and $C_{2}(m)$ which do not depend on $n$, the R.H.S. of $(5)$ is bounded above by \\ \\
$\epsilon ^{-2}\ n^{-1} \ \mathbb{E}\left(h^{2}(X_{1},\ldots,X_{m}) I_{(|h|\leq n^{\frac{3}{2}})}\right)
\left(C_{1}(m)\  + C_{2}(m) \ \log(n) \right)$
$$\leq  \ \epsilon ^{-2}\ C_{1}(m)  \ n^{\frac{-1}{3}}\ \mathbb{E}|h(X_{1},\dots,X_{m})|^{\frac{4}{3}}\qquad \qquad \qquad \qquad \qquad \qquad \qquad \qquad \qquad \qquad$$
$$\ \ +\epsilon ^{-2}\ C_{1}(m)\ \mathbb{E}\left(\ |h(X_{1},\ldots,X_{m})|^{\frac{4}{3}}\ I_{(n<|h|\leq n^{\frac{3}{2}})}\right)\ \qquad \qquad \qquad \qquad \qquad \qquad \ \ \ $$
$$\ \ +\epsilon^{-2}\ C_{2}(m)   \ n^{\frac{-1}{3}} \ \log(n)\ \mathbb{E}|h(X_{1},\dots,X_{m})|^{\frac{4}{3}}\qquad \qquad \qquad \qquad \qquad \qquad \qquad \ \ \ \ \ \ \ \ $$
$$\ \ +\epsilon^{-2}\ C_{2}(m) \ \mathbb{E}\left(|h(X_{1},\dots,X_{m})|^{\frac{4}{3}}\log|h(X_{1},\ldots,X_{m})|\ I_{(n<|h|\leq n^{\frac{3}{2}})}\right)\qquad \qquad\ \ \ $$
$$\leq  \ \epsilon ^{-2}\ C_{1}(m)  \ n^{\frac{-1}{3}}\ \mathbb{E}|h(X_{1},\dots,X_{m})|^{\frac{4}{3}}\qquad \qquad \qquad \qquad \qquad \qquad \qquad \qquad \qquad \qquad$$
$$\ \ +\epsilon ^{-2}\ C_{1}(m)\ \mathbb{E}\left(\ |h(X_{1},\ldots,X_{m})|^{\frac{4}{3}}\ I_{(|h|> n)}\right)\qquad \qquad \qquad \qquad \qquad \qquad \qquad \ \ \ \ $$
$$\ \ + \epsilon^{-2}\ C_{2}(m)   \ n^{\frac{-1}{3}} \ \log(n)\ \mathbb{E}|h(X_{1},\dots,X_{m})|^{\frac{4}{3}}\qquad \qquad \qquad \qquad \qquad \qquad \qquad \qquad \  $$
$$\ \ +\epsilon^{-2}\ C_{2}(m) \ \mathbb{E}\left(|h(X_{1},\dots,X_{m})|^{\frac{4}{3}}\log|h(X_{1},\ldots,X_{m})|\ I_{(|h|>n)}\right)\qquad \qquad\ \ \ \qquad $$
$$\longrightarrow 0, \ \ \ \ \ \textrm{as} \ n\rightarrow \infty.\qquad \qquad \qquad \qquad \qquad \qquad \qquad \qquad \qquad \qquad \qquad \qquad \qquad    $$
\\
Thus $J_{3}(n)=o_{P}(1)$. This also completes the proof of (4), and hence also that of  Theorem 3. Now, as already noted above, the proof of Theorems 1 and 2 follow from Theorem 3 and Lemmas 1 and 2.\\
\\
\textbf{Remark 5}. Studying a $U$-statistics type process that can be written as a sum of three $U$-statistics of order $m=2$, Cs\"{o}rg\H{o},   Szyszkowicz and Wang in \cite{cs1} proved that   under the slightly more relaxed  condition that $\mathbb{E}|h(X_{1},\ldots,X_{m})|^{\frac{4}{3}}<\infty$, as $n\rightarrow \infty$, we have\\
$$n^{\frac{-3}{2}}\max_{1\leq k\leq n}\sum_{1\leq i<j\leq k}\textbf{(}h(X_{i},X_{j})-\tilde{h}_{1}(X_{i})-\tilde{h}_{1}(X_{j})\textbf{)}=o_{P}(1).$$
In the proof of the latter, the well known  Doob maximal inequality for martingales was used, which gives us a sharper bound. The just mentioned inequality is not  applicable for the processes in Theorems 1 and 2, even for  $U$-statistics of order 2. The reason for this is that the inside parts of the absolute values of  $J_{s}(n),\ s=1,2,3,$ are not martingales. Also, since  $\sum _{C(k,m)} \textbf{(} h(X_{i_{1}}, \ldots,X_{i_{m}})-\tilde{h}_{1}(X_{i_{1}})- \ldots-\tilde{h}_{1}(X_{i_{m}})\textbf{)}$, for $m>2$, no longer form a martingale, it seems that  the Doob maximal inequality is not applicable  for the process $$n^{-m+\frac{1}{2}}\ \max_{1\leq k\leq n}\sum _{C(k,m)} \textbf{(} h(X_{i_{1}}, \ldots ,X_{i_{m}})-\tilde{h}_{1}(X_{i_{1}})- \ldots -\tilde{h}_{1}(X_{i_{m}})\textbf{)},$$which is an extension of the $U$-statistics parts of the process used by   Cs\"{o}rg\H{o}, Szyszkowicz and Wang in \cite{cs1} for $m=2$.
\\
\\
\par
Due to the nonexistence of the second moment of the kernel of the underlying $U$-statistic in the following example, the weak convergence result of Theorem A fails to apply. However, using Theorem 1 for example, one  can still derive  weak convergence results for the underlying $U$-statistic.\\ \\
\textbf{Example.} Let $X_{1},X_{2}, \dots$, be a sequence of i.i.d. random variables with the density function
\\
$$ f(x)=\left\{
          \begin{array}{ll}
            |x-a|^{-3}, & \hbox{$|x-a|\geq1,\ a\neq 0,$} \\
            0\ \ \ \ \ \ \ \ \ \   , & \hbox{elsewhere.}
          \end{array}
        \right.$$
Consider the parameter $\theta=\mathbb{E}^{m}(X_{1})=a^{m}$, where $m\geq 1$ is a positive integer, and the kernel $h(X_{1},\ldots,X_{m})=\prod_{i=1}^{m} X_{i}$.
Then with $m, n$ satisfying $n\geq m $, the corresponding U-statistic is
$$U_{n}={n \choose m}^{-1}
                                                                                                  \sum_{C(n,m)} \prod_{j=1}^{m} X_{i_{j}}.$$
Simple calculation shows that $\tilde{h}_{1}(X_{1})=X_{1}\ a^{m-1}\ -\ a^{m}$.\\
\\
It is easy to check that $\mathbb{E}\left(|h(X_{1},\dots,X_{m})|^{\frac{4}{3}}\log|h(X_{1},\ldots,X_{m})|\right)<\infty$ and that $\tilde{h}_{1}(X_{1})\in DAN$ (cf. Gut \cite{gut}, page 439). In order to apply Theorem 1 for this $U$-statistic, define\\
\\ $$U_{[nt]}^{*}=\left\{
                    \begin{array}{ll}
                      \ \  \ 0 \qquad \qquad  \ \qquad \qquad \qquad \qquad \qquad \qquad \ \ \  ,& \hbox{$0\leq t <\displaystyle{\frac{m}{n}}, $}\\ \\
                     {\frac {{ {[nt] \choose m}^{-1}\sum_{C([nt],m)} \prod_{j=1}^{m} X_{i_{j}}\ -\ a^{m}}}{  \textbf{(}\sum_{i=1}^{n}(X_{i}
                                                                                                  \ a^{m-1}\ -\ a^{m})^{2}\textbf{)}^{\frac{1}{2}}}} \  \ \ \ \ \ \ \ \ \ \ \qquad\ , & \hbox{$\displaystyle{\frac{m}{n}}\leq t\leq 1.$}\\
                    \end{array}
                  \right.$$\\ \\
Then, based on (c) of Theorem 1, as $n\rightarrow\infty$, we have
$$\frac{[nt]}{m} \ U_{[nt]}^{*} \longrightarrow_{d} W(t) \ \ \textrm{on} \ (D[0,1],\rho),$$
where $\rho$ is the sup-norm metric for functions in $D[0,1]$ and $\{W(t),\ 0\leq t\leq 1\}$ is a standard Wiener process. Taking $t=1$ gives us a central limit theorem for the pseudo-self-normalized $U$-statistic $$U^{*}_{n}=\frac{{n \choose m}^{-1}
                                                                                                  \sum_{C(n,m)} \prod_{j=1}^{m} X_{i_{j}}-a^{m}}{\textbf{(}\sum_{i=1}^{n}(X_{i}\
                                                                                                  a^{m-1}\ -\ a^{m})^{2}\textbf{)}^{\frac{1}{2}}}.$$
i.e., as $n\rightarrow\infty$,  we have\\
$$\frac{n}{m} \ U^{*}_{n}\longrightarrow_{d} N(0,1). $$\\ \\ \\ \\
\textbf{Acknowledgments.} The author wishes to thank Mikl\'{o}s Cs\"{o}rg\H{o}, Barbara Szyszkowicz and Qiying Wang for calling his attention to a preliminary version of their paper \cite{cs1} that inspired the truncation arguments of the present exposition. This work constitutes a part of the author's Ph.D. thesis in preparation, written under the supervision and guidance of Mikl\'{o}s Cs\"{o}rg\H{o} and Majid Mojirsheibani. My special thanks to them for also reading preliminary versions of this article, and for their instructive comments and suggestions that have much improved the construction and presentation of the results .


\newpage
\begin{thebibliography}{9999}
\bibitem{boro}Borovskikh, Yu. V. (1996). \emph{U-statistics in Banach Spaces.} VSP,
Utrecht.
\bibitem{ana} Cs\"{o}rg\H{o}, M., Szyszkowicz, B. and Wang, Q. (2003). Donsker's theorem for self-normalized parial sums processes. \emph{The Annals of Probability} \textbf{31}, 1228-1240.
\bibitem{cs} Cs\"{o}rg\H{o}, M., Szyszkowicz, B. and Wang, Q. (2004). On Weighted
Approximations and Strong Limit Theorems for Self-normalized
Partial Sums Processes. In\emph{ Asymptotic methods in Stochastics},
489-521, Fields Inst. Commun.\textbf{44}, Amer. Math. Soc., Providence, RI.
\bibitem{cs5}Cs\"{o}rg\H{o}, M., Szyszkowicz, B. and Wang, Q. (2008). On weighted approximations in $D[0,1]$ with application to self-normalized partial sum processes.
\emph{Acta Mathematica Hungarica} \textbf{121 (4)}, 307-332.
\bibitem{cs1} Cs\"{o}rg\H{o}, M., Szyszkowicz, B. and Wang, Q. (2008). Asymptotics of studentized U-type processes for changepoint problems. \emph{Acta Mathematica Hungarica} \textbf{121 (4)}, 333-357.
\bibitem{Gin} Gin\'{e}, E. , G\"{o}tze, F. and Mason D. M. (1997). When is the student t-statistic asymptotically Normal? \emph{The Annals of Probability} \textbf{25}, 1514-1531.
\bibitem{gut} Gut, A. (2005). \emph{Probability: A Graduate Course.} Springer.
\bibitem{hof} Hoeffding, W. (1948). A class of statistics with asymptotically normal distribution. \emph{Ann. Math. Statist.} \textbf{19}, 293-325.
\bibitem{mil} Miller, R. G. Jr. and Sen, P. K. (1972). Weak convergence of U-statistics and Von Mises' differentiable statistical functions. \emph{Ann. Math. Statist.} \textbf{43}, 31-41.
\bibitem{serf}Serfling, R. J. (1980). \emph{Approximation Theorems of Mathematical
Statistics}. Wiley, New York.
\end{thebibliography}
\end{document}